\documentclass[12pt,preprint,authoryear,round]{imsart}
\RequirePackage[colorlinks,citecolor=blue,urlcolor=blue]{hyperref}
\usepackage{amssymb,amscd,amsthm, verbatim,amsmath,color,fancyhdr, mathrsfs}
\usepackage{graphicx,tabularx,booktabs}
\usepackage{turnstile}
\usepackage{fancyhdr}
\usepackage{caption,subcaption}
\usepackage{url}
\usepackage{color,verbatim,enumerate}
\usepackage[normalem]{ulem}
\usepackage{marginnote} 
\usepackage{multirow}
\usepackage[letterpaper, left=2.5cm, right=2.5cm, top=2.5cm,
bottom=2.5cm,dvips]{geometry}
\usepackage{tikz}
\usepackage{natbib}

\usepackage[ruled, 
lined, 
commentsnumbered]{algorithm2e}

\theoremstyle{plain} 

\theoremstyle{definition}

\newcommand\independent{\protect\mathpalette{\protect\independenT}{\perp}}
\def\independenT#1#2{\mathrel{\rlap{$#1#2$}\mkern2mu{#1#2}}}
\newcommand\notindependent{\!\perp\!\!\!\!\not\perp\!}

\begin{document}
		\begin{frontmatter}
			\title{Algebraic Statistics in Practice: \\ Applications to Networks}
			\runtitle{Algebraic Statistics in Practice}
			
			\author{Marta Casanellas, Sonja Petrovi\'c, and Caroline Uhler}
			\runauthor{Casanellas, Petrovi\'c, Uhler }
			

\begin{abstract}
Algebraic statistics uses tools from algebra (especially from multilinear algebra, commutative algebra and computational algebra), geometry and combinatorics to provide insight into knotty problems in mathematical statistics.  In this survey we illustrate this on three problems related to networks, namely network models for relational data, causal structure discovery and phylogenetics. For each problem we give an overview of recent results in algebraic statistics with emphasis on the statistical achievements made possible by these tools and their practical relevance for applications to other scientific disciplines.
\end{abstract}

	\end{frontmatter}	



\section{Introduction}
Algebraic statistics is a branch of mathematical statistics that focuses on the use of algebraic, geometric and combinatorial methods in statistics. The term ``Algebraic Statistics'' itself was coined as the title of a book on the use of techniques from commutative algebra in experimental design~\citep{Pistone2001}. An early influential paper \citep{DS98} connected the problem of sampling from conditional distributions for the analysis of categorical data to commutative algebra, thereby showcasing the power of the interplay between these areas.
In the two decades that followed, growing interest in applying new algebraic tools to key problems in statistics has generated a growing literature.

The use of algebra, geometry and combinatorics in statistics did not start only two decades ago. Combinatorics and probability theory have gone hand-in-hand since their beginnings. The first standard mathematical method in statistics may be the Method of Least Squares, which has been used extensively since shortly after 1800 and relies heavily on systems of linear equations. Non-linear algebra has played a major role in statistics since the 1940s; see for example~\citet{Wilks_1946}, \citet{Votaw_1948}, \citet{James_1954}, \citet{Andersson_1975}, \citet{Bailey_1981}, and \citet{Jensen_1988}. In addition, the development of the theory of exponential families relied heavily on convex geometry~\citep{Barndorff_Nielsen}. However, with ~\citet{DS98} and \citet{Pistone2001}, new algebraic disciplines including modern computational algebraic geometry and commutative algebra were introduced in statistics. In this review, we concentrate on the developments in algebraic statistics in the last two decades and in particular on applications to networks.

The analysis of networks  as relational data and to represent probabilistic interactions between variables is becoming increasingly popular, with applications in fields including the social sciences, genomics, neuroscience, economics, linguistics and medicine.
Theoretical and algorithmic developments for exploring such datasets
 are found at the intersection of statistics, applied mathematics,  and machine learning.
 In this review we focus on some of the key statistical problems and their solutions using  algebraic techniques 
  in three application areas:
 network models (based on relational data encoded as \emph{observations on the edges} of a random network), causal structure discovery (based on multivariate data encoded as \emph{observations on the nodes} of an unknown underlying causal network), and phylogenetics (a particular network structure discovery problem where the underlying network is a tree with latent variables).

 Section \ref{chp:sonja}
focuses on statistical models for relational data,
typical  uses of which arise in the social and biological sciences. In these applications, nodes in the network may represent individuals, organizations, proteins, neurons, or brain regions, while  links 
 represent observed relationships between the nodes, such as personal or organizational affinities, social/financial relationships, binding between proteins or physical links between brain regions.
A key problem  in this area is to test whether a proposed statistical model fits the data at hand; such a test 
typically involves generating a sufficiently large and generic sample of networks from the model and comparing it to the observed network. Perhaps somewhat surprisingly, algorithms for sampling networks with given network statistics for goodness-of-fit testing are often efficiently encoded by algebraic constraints.  In Section~\ref{chp:sonja}, we outline how techniques from commutative algebra and combinatorics 
 are applied to 
 this problem
   for several families of network models for which a formal test is otherwise unavailable.

In Section~\ref{chp_3}, we turn to applications where the network structure cannot directly be observed and we only have access to observations on the nodes of the network. Such applications range from data on consumer behavior to click statistics for ads or websites, DNA sequences of related species, gene expression data, etc. The use of such data to gain insight into complex phenomena requires characterizing the relationships among the observed variables. Probabilistic graphical models explicitly capture the statistical relationships between the variables as a network. A good representation of a complex system should not only enable predicting the state of one component given others, but also the effect that local operations have on the global system. This requires causal modeling and making use of interventional data. In Section~\ref{chp_3}, we discuss the role that algebraic and discrete geometry play in analyzing prominent algorithms for causal structure discovery and in developing the first provably consistent algorithms for causal inference from a mix of observational and interventional data.

In Section~\ref{chp:marta}, we discuss a particular directed network model, namely phylogenetic trees, for evolutionary reconstruction. Algebra and related areas have always been present in the study of evolutionary processes, but have played minor roles relative to combinatorics or optimization. However, since the beginning of this century, the developments in algebraic statistics have given rise to techniques with a major impact on three different problems in phylogenetics: model selection, model identifiability, and phylogenetic reconstruction. 
  Models 
  that best fit the data 
 should only be selected among those whose parameters are identifiable and hence understanding model identifiability is crucial. The final step given an evolutionary model and data is to reconstruct the phylogenetic tree and infer the evolutionary parameters. In Section~\ref{chp:marta}, we explain how algebraic techniques can be used to address these problems and discuss their applications to complex evolutionary models and phylogenetic networks.

\section{Network models for relational data} 
\label{chp:sonja}

\begin{figure}
	\includegraphics[scale=0.3]{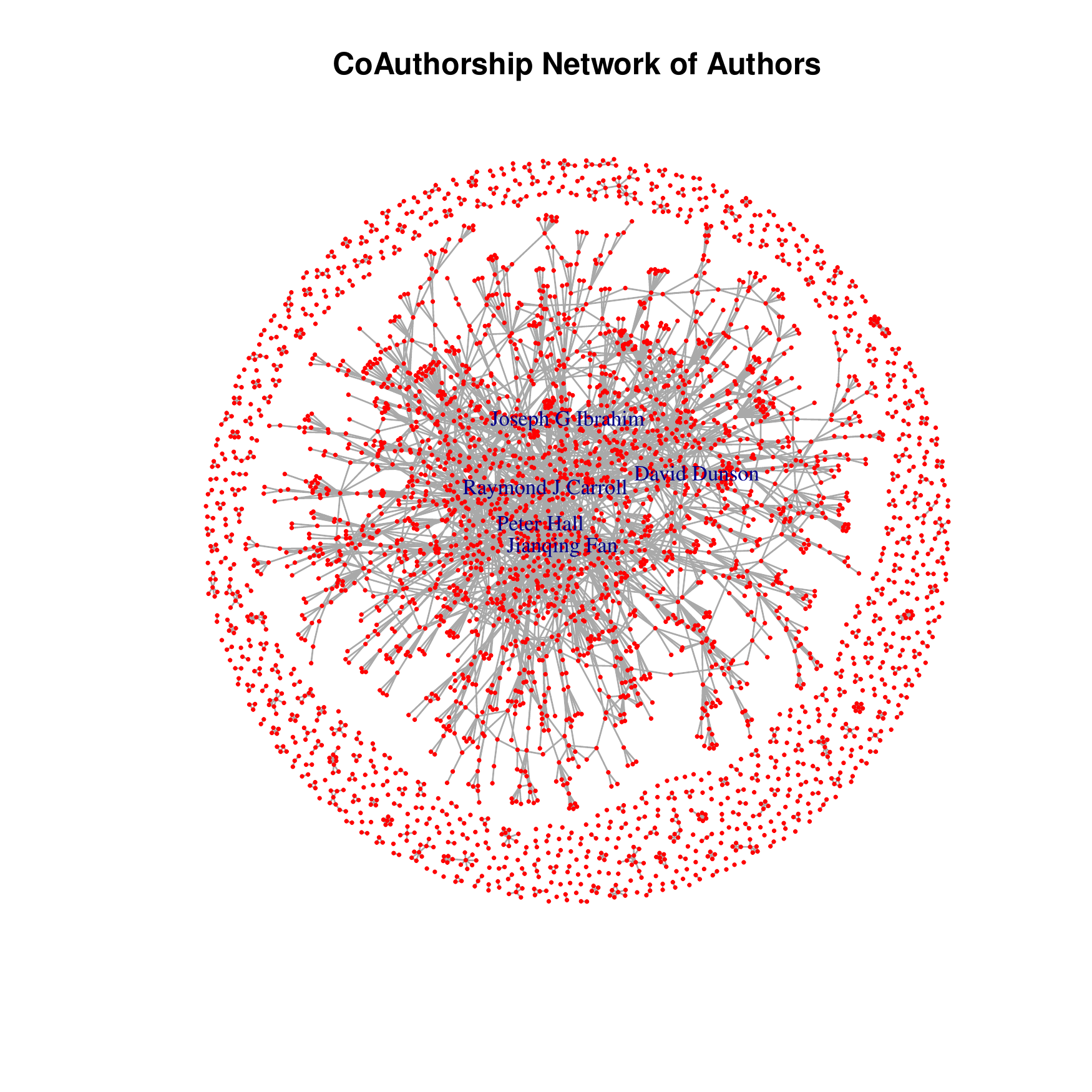}
	\includegraphics[scale=0.3]{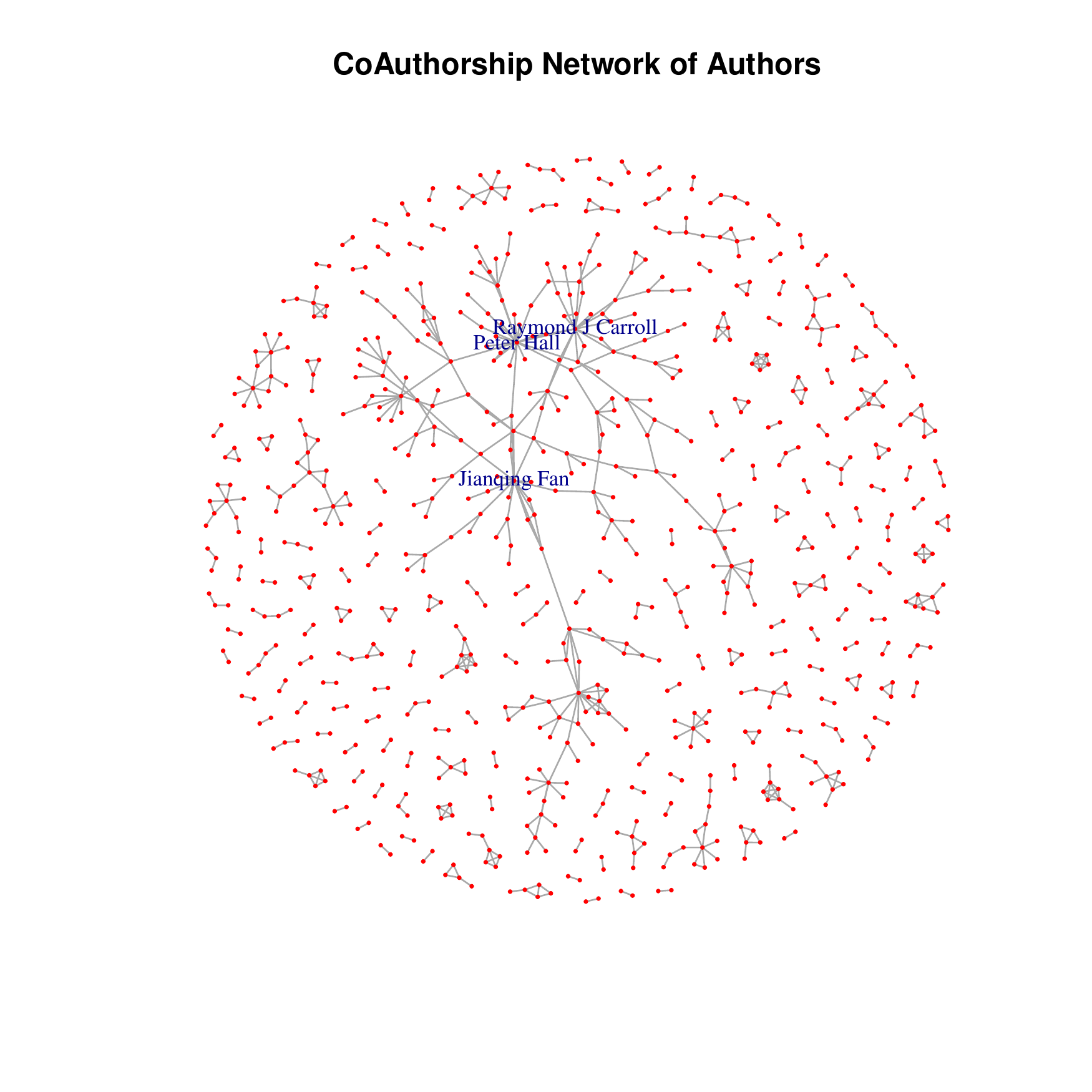}
	\includegraphics[scale=0.3]{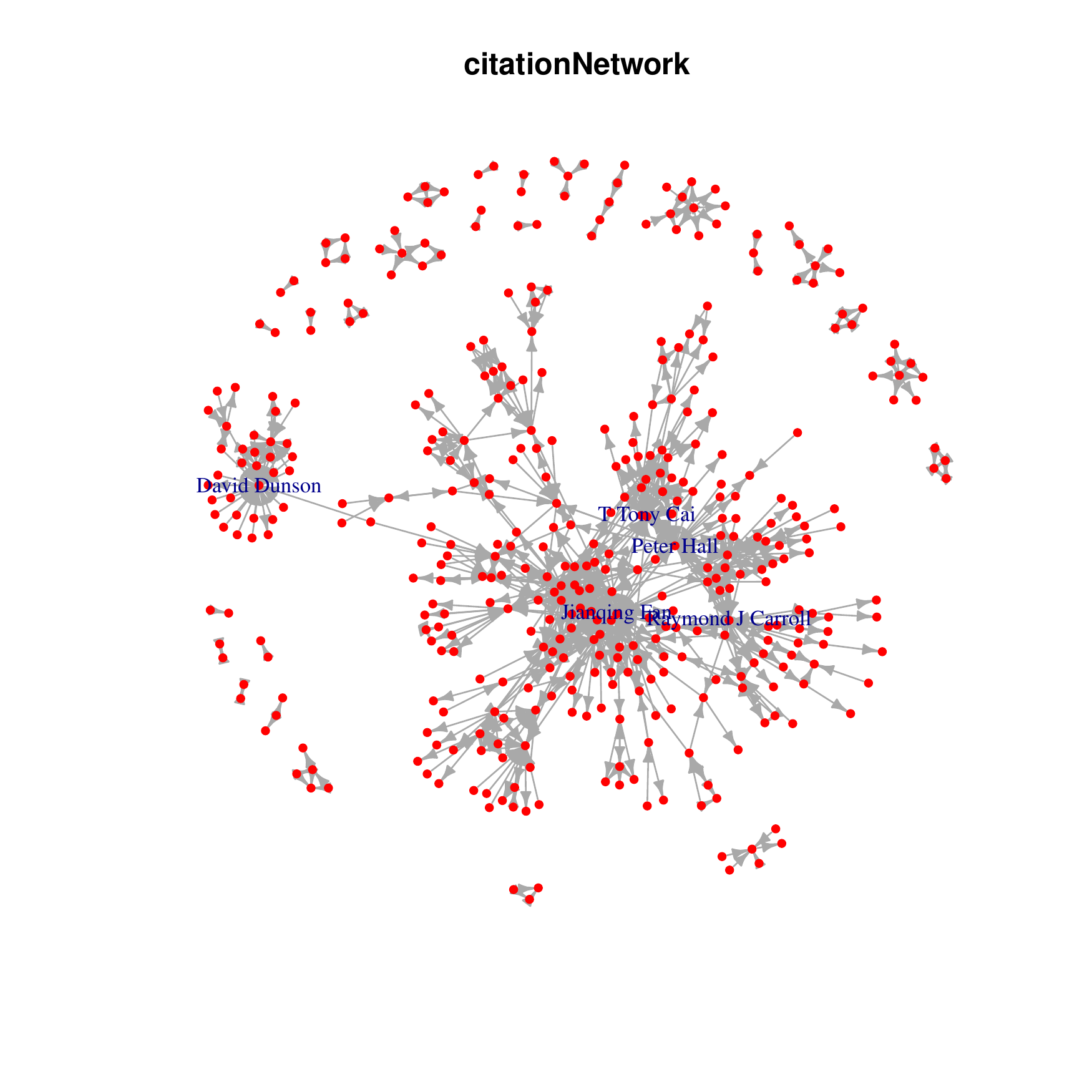}
	\vspace{-0.8cm}
	\caption{Author networks constructed from the data collected by Ji and Jin \citep{JiJinAuthorsAOAS}. In  ``Coauthorship network A'' (left) there is an undirected edge between nodes $i$ and $j$ if authors $i$ and $j$ coauthored at least $2$ papers. In  ``Coauthorship network B'' (center) there is an undirected edge between nodes $i$ and $j$ if authors $i$ and $j$ coauthored at least $1$ paper. In  ``Citation network" (right), there is a directed edge from author $i$ to author $j$ if  $i$ cited  at least $1$ paper by  $j$.
}\label{fig:authorsData}
\end{figure}

Network models for relational data, that is, various types of interactions between a fixed set of entities, such as neurons, proteins, people, or corporations, have grown in popularity in  recent decades. The interactions can be directed (e.g.,~affinity or one-way influence) or undirected (e.g.,~mutual affiliation), and may  be counted with multiplicity or weight. 
Consider  two recently-collected data sets on statisticians who publish in five top-rated journals~\citep{JiJinAuthorsAOAS}. 
The data can be represented  as a bipartite graph of authors and papers, in which a link exists between nodes $i$ and $j$ if author $i$ wrote paper $j$, or as a citation network, in which a directed edge from $i$ to $j$ denotes that paper $i$ cites paper $j$, or as collapsed coauthorship or citation networks among authors; see Figure~\ref{fig:authorsData}. 

Taking a model-based approach, we study the effects of various types of author interactions on network analysis and inference by concentrating on goodness of fit of a network model. 
  This  is central for 
  estimating network features, appropriately simulating data,  and correctly interpreting the results. In addition, it is considered to be very challenging in the network science community due to both the size of the networks in many applications as well as their sparsity and particular structure.  
Methods rooted in algebraic statistics help answer such questions efficiently and reliably  for a variety of network models.

Relational data can be modeled using random graphs, in which the interactions are modeled as random variables. This results in a statistical model with random directed or undirected edges on a fixed set of nodes. 
 There is a rich literature on different random graph models, starting from the classical Erd\"{o}s--R\'{e}nyi graphs  \citep{erdos1961evolution}, exponential random graph models \citep{HL81}, and Markov graphs  \citep{frank1986markov}, to  models that capture more intricate relational behavior, such as stochastic blockmodels \citep{holland1983stochastic},  latent space models \citep{hoff2002latent}, and mixed membership stochastic blockmodels \citep{airoldi2009mixed}; see also   \citep{goldenberg2010survey}. 
The question  whether any of these  models provides an adequate fit to  data  has received relatively little attention.  

Here we consider 
 the broad and flexible class of \emph{exponential family models for random graphs}, also known as {ERGMs}. 
 To specify an ERGM, one first selects a vector of network characteristics $T(g) \subset \mathbb{R}^p$ that represent an interpretable and meaningful summary of the network, such as the number of neighbors of each node, block membership, etc. The resulting model is the collection of probability measures $\mathcal{M}= \{p_\theta : \theta \in \Theta\}$,  indexed by points in $\Theta \subset \mathbb{R}^p$ such that for any $\theta \in \Theta$, the probability of observing a given network $G = g$ takes  the exponential form  
 \[
 	p_\theta (g) = \exp\{\langle T (g) , \theta \rangle - \psi(\theta)\},
\]
 where $\psi(\theta) = \log \sum_g \exp \{\langle T (g) , \theta \rangle \}$ is the normalizing function (also known as  the log-partition function) and $T(g)$ is  the vector of minimal sufficient statistics for $\mathcal{M}$.
 Tools from commutative algebra can be applied to construct a finite-sample test for goodness of fit of such a model to the observed network, while graph-theoretic and combinatorial considerations can render the resulting algorithms scalable and applicable in practice for large networks.

\subsection{Testing  model fit: state-of-the-art} 
\label{sec:issues}

Studies devoted to goodness-of-fit tests for network models  fall into  two categories. 

\emph{Heuristic tests} are based on graphical comparisons between observed statistics and the corresponding statistics obtained from the fitted model, see \citep{H:03}, \citep{NtwkApproxFittingDynamics15}, \citep{hunter2008goodness}. 
Given an observed graph $g_{obs}$, the goal is to evaluate how well a  model $\mathbb P_{\theta}(G)$  fits $g_{obs}$.
 Let $s(g)$ be a vector of network statistics; most popular ones include  the number of edges, triangles, or two-stars in $g$,   the vector of counts of neighbors of every node (the degree sequence), or other summaries of a node's connectedness or centrality in $g$. 
 The graphical method proceeds by computing a maximum likelihood estimator (MLE) $\hat \theta$ of $\theta$ and simulating several graphs $g^{1},\ldots, g^{B}$ from $\mathbb P_{\hat \theta}$. Departures from the model are detected by comparing the sample distribution of $s(g^{1}), \ldots, s(g^{B})$ with the observed value $s(g_{obs})$.   
Central to the graphical method is the choice of complementary statistics $s(g)$ used for evaluating the fit.  
While widely used,  graphical tests  
 have two limitations: First, they are not based on any formal discrepancy measure between the model and observed network, since the choice of $s(g)$ is arbitrary. Second, the distribution of the complementary statistics is unknown under the null hypothesis, so calibration and formal Type I error rates are difficult to obtain. 

\emph{Asymptotic tests} are a natural alternative and rely on   formal testing criteria for evaluating  model fit. However, classical test criteria such as the log-likelihood ratio, AIC, or BIC, cannot be directly applied to general network models, mainly because the usual asymptotics do not apply to models other than very simplistic ones. This is due to the fact that the iid assumption on the random edges does not hold,  which dismantles results on asymptotic distributions of various test statistics. In addition, in many network models the number of parameters increases with the number of nodes. This issue was first pointed out in  \citep{FW81} and noted also in several later works \citep{KrivitskyKolaczyk15, hunter2008goodness, HL81,p1asymptotics, NtwkApproxFittingDynamics15, CDS11}. In addition, many commonly-used ERGMs  suffer from the lack of a natural notion of projectability \citep{shalizi2013consistency}, which relates the marginal distribution of a network on $p$ nodes to the same model on $p+1$ nodes, essentially 
 ruling out  consistency of MLEs. 

To remedy these issues, one can derive modified asymptotic distributions, when they exist, of various test statistics for special cases. For example, \citet{yan2014model} consider testing the degree-corrected blockmodel with the usual stochastic blockmodel; 
\citet{wang2017likelihood} derive an asymptotic Gaussian distribution of the likelihood ratio test statistic for selecting between two stochastic blockmodels with different number of communities; 
\citet{gao2017testing} consider  testing an Erd\"{o}s-R\'{e}nyi model against a stochastic blockmodel,  construct a chi-square-like test statistic using a combination of edge, 2-star, and triangle counts, and show that its limiting distribution is a chi-square distribution; \citet{lei2016goodness} constructs a goodness-of-fit test for the stochastic blockmodel by using the extreme eigenvalues of a certain residual matrix as a test statistic and deriving its asymptotic distribution; similarly, \citet{banerjee2017optimal} derive a central limit theorem for linear spectral statistics for testing an Erd\"{o}s-R\'{e}nyi model against a two-block blockmodel.
A common limitation of these studies is that the asymptotic distributions are derived in specialized asymptotic regimes that may not hold in practice and are difficult to verify, given a single sample of a network. For instance, in \citep{lei2016goodness} the asymptotic null distribution of the test statistic requires that the entries of the estimated edge probabilities be uniformly bounded away from $0$ and $1$, which rules out certain types of sparse networks.

\subsection{From networks to contingency tables: log-linear models}

 Network data on $p$ nodes can be naturally summarized by a contingency table of format $p\times p\times i_1\times\dots\times i_k$,  classifying the type of a relationship (directed, undirected, block-dependent, etc.) that holds for each dyad in the graph. 
This representation means that certain ERGMs can be represented by equivalent models for contingency tables that have a long history in the statistics literature. 
The models amenable to such a representation are called \emph{log-linear ERGMs}, their   vector of sufficient statistics is a linear function of the network. 
For such models there exists a matrix $A$ such that  $T(g) = Ag$, where the network $g$ has been flattened to vector format. 
Log-linear ERGMs  encompass many of the popular models in use today, including all undirected and directed degree-based models (e.g., the $\beta$-model \citep{CDS11,RPF:11}), stochastic blockmodels or SBMs with or without mixed membership (but with known block assignment  \citep{HL81,FienbergMeyerWasserman1985block,airoldi2009mixed}),   combinations of these (e.g., the degree-corrected SBM \citep{karrer2011stochastic}), and extensions of any of these models using covariates \citep{YJFL-InferenceDirectedNtwkCovariates-JASA}. 

The  connection to 
contingency tables 
 dates back to three seminal papers from the 1980's, namely \citep{FienbergWasserman1981categorical,FienbergMeyerWasserman1985block,FW81}, which consider some  (very novel at the time and still very popular today) models for relational data. 
  By viewing the network representation of the data as a union of independent dyads that can appear in various configurations, they express in table format  a set of models that are now considered canonical  under the ERGM framework. 
The first advantage of this viewpoint,  also pointed out in these early works, is that the MLE can efficiently and accurately be computed using iterative proportional fitting,  thus avoiding the usual convergence issues that are the main drawback of MCMC approaches typically used for ERGMs; see e.g.~\citep{ergm}. 
The second advantage became apparent in the 2000s with the development of tools from algebraic statistics for contingency tables: a  generating set of a polynomial ideal can be translated to a set of networks that preserve an ERGM's sufficient statistics, then used as input to a sampling algorithm that provides a reference set for testing model fit.  Coupled with a valid discrepancy measure for model fit also ported from the contingency table literature into networks, this approach  solves the issues outlined in \ref{sec:issues}. 

\subsection{Goodness-of-fit testing for log-linear ERGMs}

Let $\mathcal M_T$ be a log-linear ERGM, where $T$ denotes the vector of sufficient statistics. A canonical way to test model fit is to compute the exact  $p$-value conditional on the sufficient statistics for the null hypothesis that $p_{\hat{\theta}}(g)$ lies in the model $\mathcal M_T$, where $\hat{\theta}$ is the MLE, against the general alternative (see~\citep{FW81} for further motivation). The $p$-value is computed by comparing the observed network $g$ against all other networks whose sufficient statistics are the same; this set, 
$$\mathcal F_T(g) := \{g': T(g')=T(g)\}$$ is called \emph{the fiber of $g$ under the model $\mathcal M_T$}. In virtually all instances of interest for applications, the fiber is too large to enumerate, so one resorts to sampling from it.

\begin{figure}
	\includegraphics[scale=0.8]{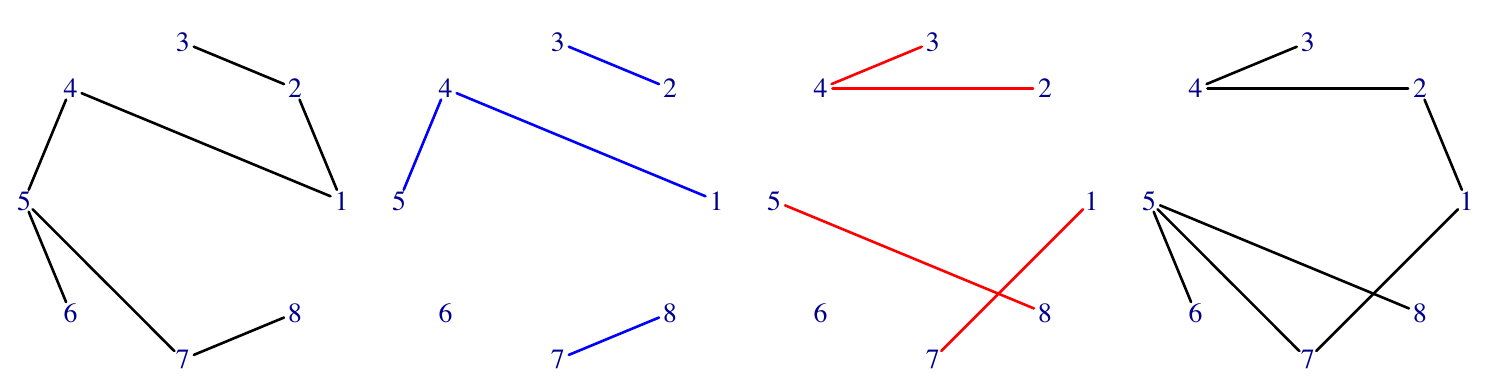}

	\includegraphics[scale=0.8]{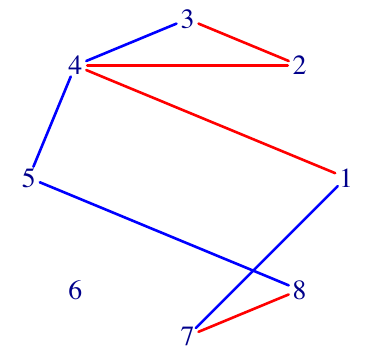}
	\caption{Example of a move for the ERGM with $T(g)=(d_1,\dots,d_8)=(2,2,1,2,3,1,2,1)$, where $d_i$ is the number of neighbors of node $i$.  
	{\bf Top}, in order: a starting graph $g$, set of blue edges to be removed from $g$, set of red edges to be added to $g$, and the resulting graph $h$. 
	{\bf Bottom}: the move $b$ represented as a bicolored graph, with blue edges carrying weight $-1$ and red $+1$. 
	Since blue edges contribute `negative' neighbors,  $T(b)=0$ and thus $T(h)=T(g)$. 
}\label{fig:moves}
\end{figure}

To sample from this conditional distribution for any log-linear model, \citet{DS98} introduce a notion of a basis that can be used as input to the Metropolis-Hastings algorithm. 
 In the context of networks,  
 a \emph{Markov basis} of the log-linear ERGM $\mathcal M_T$ is any set of networks $\mathcal B= \{b_1,\dots,b_n\}$ 
 for which 
$
	T(b_i) =0 
$ 
and such that 
for any given network  $g$ 
and any  $h\in\mathcal F_T(g)$, there exist $b_{i_1},\dots,b_{i_N}\in \mathcal B$ that can be used to reach $h$ from $g$, i.e.,  
\[
	g + b_{i_1} +\ldots+ b_{i_N} = h,
\]
while walking through elements of the fiber, meaning that each partial sum $u+\sum_{j=0}^N b_{i_j}$, for any $j=1,\dots,N$, represents a valid network; see Figure~\ref{fig:moves}. 
Note that $T(u)=T(u+b_i)$ means that adding a move $b_i$ to any network  does not change the values of the sufficient statistics, so to remain in the fiber,  we  only need to ensure  that adding a move did not produce negative entries in the vector, as the count of edges in a graph cannot be negative. 
The resulting Markov chain is irreducible, symmetric, and aperiodic;  \citet[Algorithm 1.13]{DSS09} outlines a vanilla implementation. 

The connection to commutative algebra translates each move $b_i$ into a binomial: a difference of products of indeterminates, each corresponding to a cell in the contingency table.
   In the example from Figure~\ref{fig:moves}, the depicted move can be written as $e_{17}e_{24}e_{34}e_{58}-e_{14}e_{23}e_{45}e_{78}$, where  $e_{ij}$ is the indeterminate representing the dyad $\{i,j\}$. The move is thus a polynomial in the random dyads. 
 This translation is straightforward, but leads to a fundamental and surprising result: a set of moves is a Markov basis if and only if the corresponding binomials generate the toric ideal defined by $T$ \citep{DS98}.  Consequently, each log-linear model \emph{has} a finite Markov basis, by the Hilbert basis theorem from algebra; and all the basis elements can be computed, using combinatorial tools for computing bases of toric ideals. 

Markov bases are a popular theoretical construct in algebraic statistics, but in practice  pose  serious challenges in particular pertinent to large networks, and in general for large sparse contingency tables. 
One  is that they are complicated to compute a priori and that algebra produces many moves inapplicable to the  observed data. To circumvent this difficulty, \citet{GPS16,karwa2016exact,GPS19} implement a \emph{dynamic} algorithm for generating Markov elements for the $\beta$ and $p_1$ models, some of the basic variants of the stochastic blockmodels, and combinations of these, and embed them into a Metropolis-Hastings algorithm to provide a scalable exact conditional test for model fit. 
Another concerns the mixing time of the Markov chain constructed using Markov bases, as any Markov chain that is slow to mix will not be scalable to large networks in practice; to this end, we will only mention that there is a large body of literature in discrete mathematics that implies rapid mixing of this chain for almost all fibers, for details, see \citep{DillonMS16}. 

\vspace{0.3cm}
\noindent\textbf{Example:} Considering the  largest connected component of the citation network of authors, 
 \cite{JiJinAuthorsAOAS} count the neighbors in this directed graph 
  and propose using them for author rankings. 
 We perform an exact test of model fit for the ERGM whose sufficient statistics is this vector of neighbor counts, namely, the   $p_1$ model with dyad-dependent reciprocation. 
The test is done by running the Markov chain  described in \citet{GPS16}. After  $N=100,000$ steps, the estimated $p$-value is $0.0072$.
 As a measure of discrepancy between the observed graph and the MLE, we use the chi-square statistic. The $p$-value reported is the proportion of the sampled networks in the fiber whose chi-square value is at least as large as that of $g_{obs}$. 
\begin{figure}
	\includegraphics[scale=0.5]{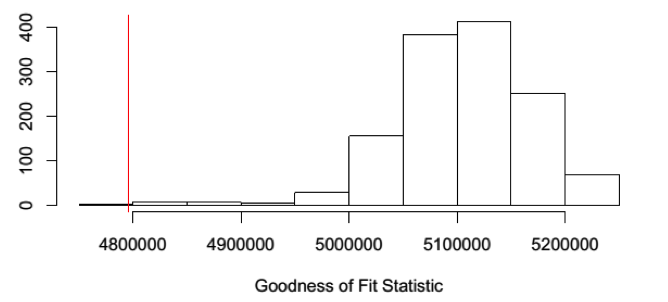}
	\qquad 
	\includegraphics[scale=0.5]{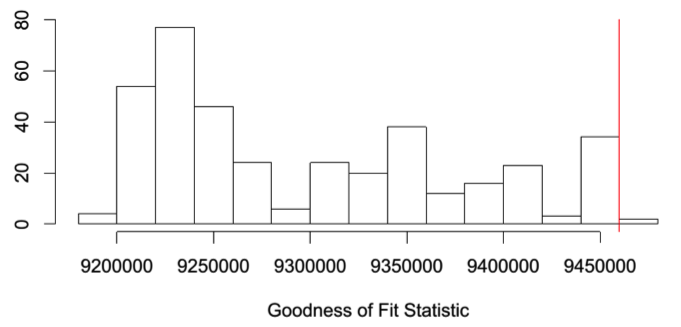}
	\caption{Testing model fit for the ERGMs with node degrees a sufficient statistics. Left: histogram of the chi-square statistics, as a measure of discrepancy, for 100,000 networks in the fiber of the coauthorship network. Right: same information for the citation network. The vertical line indicates the observed value of the chi-square statistic. 
}\label{fig:authorsDataGoF}
\end{figure}
 This result indicates that the $p_1$ model does not fit the citation network of authors, and therefore the network may posses transitive effects and the dyads may  not be independent. Similarly, we perform an exact test of model  fit for the $\beta$ model in the largest connected component of the coauthorship network $A$. The $p$-value from the goodness-of-fit test obtained by running   the Markov chain  for $N=100000$ steps is $0.997$. 
 This suggests that node degrees are a (almost surprisingly) good  summary of the graph and that the degree of an author could be used to determine an author ranking. 
 However, this graph was obtained by \emph{thresholding} the original data (a popular technique in network analysis used to avoid multiple edges), as well as by \emph{reducing} multi-author papers to pairs of authors.
These two tests are summarized in Figure~\ref{fig:authorsDataGoF}, which depicts the sampling distributions of the chi-square statistic. 
That the Markov chains converged fairly well was checked using the usual MCMC diagnostics in {\tt R}. 

\subsection{Generalizations to weighted graphs} 
The previous example opens up several interesting questions: how can we preserve the underlying data structure and still use an interesting network model with scalable estimation and goodness-of-fit methods?  
  \citet{KarwaPetrovic:AOASauthors} argue that thresholding and reducing to a graph is not necessary; one can instead work with a hypergraph representation of the data, which preserves more of the coauthorship structure than the network representation. 
For $I$ authors, $J$ research areas and $K$ journals, consider an $I \times I \times J \times K$ contingency table whose $(i, i', j, k)$ entry counts the number of times author $i$ cites author $i'$ in research area $j$ and journal $k$. A similar representation can be obtained for the coauthorship network, where  we count the number of times authors $i$ and  $j$ wrote a joint paper. These representations preserve the citation and coauthorship count data.  
We can then collapse the table to an $I \times I$ author-by-author table and fit log-linear models to the citation counts. 
In essence, we seek to avoid thresholding, as in  the generalized $\beta$ model  discussed in \citet{RPF:11} for weighted networks represented in table form. 
Generalizing to weighted or multiple graphs is straightforward in the contingency table setting, with MLE algorithms unaffected, and Markov basis algorithms becoming--perhaps surprisingly--more efficient and easier to implement. 
This opens up several lines of research on generalizing these models and enriches the network science literature with  goodness-of-fit tests for many popular ERGMs. 

By definition, log-linearity means  that the sufficient statistics are a linear function of the graph, which in turn implies dyadic independence. The assumption of dyadic independence may seem restrictive; but \citet{YJFL-InferenceDirectedNtwkCovariates-JASA} show that it includes many popular models and avoids the  degeneracy that plagues other ERGMs. In addition, \citep{karwa2016exact} develop goodness-of-fit testing methods combining the Bayesian and algebraic approaches for  mixture models of log-linear ERGMs, which do not assume dyadic independence. 


\section{Causal structure discovery}
\label{chp_3}
From random graph models, where each edge of the network is associated with a random variable, we now turn to graphical models, where each node of the network is associated with a random variable. In most applications, the underlying network is unknown and needs to be learned from data on the nodes. We here consider the problem of learning the causal relationships among the nodes.

Causal inference is the basis of scientific discovery, because it asks `why?'. The gold standard for inferring causal relationships is randomized controlled trials. 
However, in many applications running such trials 
to test for a causal effect is impractical, unethical or prohibitively expensive. So there have been large efforts to develop a theory of causal inference based purely on observational data. This began with two crucial advances made independently in the 1920s. Jerzy Neyman established a formal distinction between random variables under randomization and ordinary random variables via the potential outcome notation~\citep{Neyman}. Sewall Wright independently pioneered the use of graphs to represent cause-effect relationships using structural equation models~\citep{SWright_1, SWright_2}. However, skepticism amongst statisticians resulted in the causal interpretation of structural equation models being overlooked and almost forgotten (see~\citet{Pearl_2012} for a historical account). The reemergence of causal inference from observational data in statistics began in the 1970s and led to major contributions by \citet{Pearl}, \citet{Robins_2}, \citet{Rubin_1974, Rubin_2005}, and \citet{pc_algorithm}.


While it has in general been unethical, too expensive or even impossible to perform large-scale interventional studies, the development of genome editing technologies in biological studies~\citep{CongZhang2013} as well as the explosion of interventional data in online advertisement and education represents a unique opportunity for the development of new  causal inference methodologies. It is now possible to obtain large-scale interventional datasets relatively easily. This calls for a theoretical and algorithmic framework for learning causal networks from a mix of observational and interventional data.

In this section, we showcase how methods from algebraic geometry, combinatorics, graph theory and discrete geometry have been brought to bear in the analysis and development of causal structure discovery algorithms. In Section~\ref{SEM_Markov}, we  introduce the framework of structural equation models for causal modeling and then discuss open problems in combinatorics and graph theory related to the degree of identifiability of causal effects. In Section~\ref{algorithms}, we will review a prominent causal structure discovery algorithm. This algorithm relies on the so-called faithfulness assumption, and using algebraic geometry we will show that this assumption is very restrictive and hard to satisfy in practice. 
In Section~\ref{DAG_associahedra}, we will discuss an alternative algorithm that makes critical use of discrete geometry to overcome the limitations of the faithfulness assumption and leads to the first provably consistent algorithm for causal inference from a mix of observational and interventional data. Finally, in Section~\ref{applications} we 
discuss various open problems and related literature in algebraic statistics.

\subsection{Structural equation models and Markov equivalence}
\label{SEM_Markov}

We represent a causal network by a directed graph $G=(V,E)$ consisting of vertices $V=\{1,\dots , p\}$ and directed edges $E$ representing direct causal relationships. We make the common assumption that $G$ is a directed \emph{acyclic} graph (DAG), meaning there are no directed cycles $i_0 \to i_1\to\dots \to i_m\to i_0$, since causal effects only act forward in time. In a \emph{structural equation model}~\citep{SWright_1, SWright_2}, each node $i\in V$ is associated with a random variable $X_i$ and  is a deterministic function of its parents, denoted by $\textrm{pa}(i)$, and independent noise, denoted by $\epsilon_i$. For example, a structural equation model on the 4-node DAG $1\to 2$, $2\to 3$, $3\to 4$, $1\to 4$ is given by
\begin{equation}
\label{eq_SEM}
X_1=f_1(\epsilon_1), \quad X_2=f_2(X_1,\epsilon_2), \quad X_3=f_3(X_2,\epsilon_3), \quad X_4=f_4(X_1,X_3,\epsilon_4).
\end{equation}

\emph{Gaussian linear structural equation models} are special instances of this model class, where $X_j=\sum_{i\in\textrm{pa}(j)} a_{ij}X_i +\epsilon_j$ and the noise $\epsilon=(\epsilon_1, \dots , \epsilon_p)$ follows a Gaussian distribution $\mathcal{N}(0,D)$, where the covariance matrix is diagonal. In this case, the joint distribution of $X=(X_1, \dots , X_p)$ is a Gaussian $\mathcal{N}(0,\Sigma)$, where $\Sigma^{-1}=(I-A)D^{-1}(I-A)^T$ and $A$ is the weighted adjacency matrix of $G$ containing the causal effects $a_{ij}$. While this model is of interest for its mathematical simplicity, in many applications including genomics the linear and Gaussian assumptions are often violated and it is preferable to work with the general non-parametric model in (\ref{eq_SEM}).

A structural equation model not only encodes the \emph{observational distribution}, i.e., the distribution of $X$, but also the \emph{interventional distributions}. For instance, in the example above an intervention on node $X_3$ by setting it's value to 0 would change the distribution of the nodes $X_3$ and $X_4$, but not the others, since they are not downstream of $X_2$. Such an intervention could for example be used to model a gene knockout experiment, where the expression of certain genes is set to zero~\citep{CongZhang2013}. 

A structural equation model provides a factorization of the joint distribution, which implies certain \emph{conditional independence} (\emph{CI}) \emph{relations} through the \emph{Markov property}, namely $X_i\independent X_{\textrm{nd}(i)}\mid X_{\textrm{pa}(i)},$
where $\textrm{nd}(i)$ denotes the non-descendents of node $i$; see, e.g.~\citet{Lauritzen1996} for an introduction to graphical models. A standard approach for causal structure discovery is to infer CI relations from the sample distribution and then infer the DAG from these relations. However, in general a DAG is not identifiable, since multiple DAGs can encode the same set of CI relations; such DAGs are called \emph{Markov equivalent}. \citet{Verma_Pearl} provided a graphical characterization of when two DAGs are Markov equivalent, namely when they have the same \emph{skeleton} (i.e., undirected edges) and \emph{immoralities} (i.e.~induced subgraphs of the form $i\to j \leftarrow k$). 

Since from observational data it is only possible to identify a DAG up to its Markov equivalence class (MEC), it is important to study the sizes of MECs and their distribution. However, while a recurrence relation for the number of DAGs on $p$ nodes is known~\citep{Robinson}, no such formula is known for MECs. \citet{Gillispie} enumerated all MECs up to 10 nodes; see Table~\ref{table_MECs}. The first row shows that the number of MECs grows very quickly in the number of nodes~$p$. The second row shows the ratio of the number of MECs to the number of DAGs, suggesting that this sequence converges to $\approx 1/4$. This combinatorial conjecture would have important consequences for causal inference, since it would imply that on average a Markov equivalence class consists of about 4~DAGs, meaning that in general only very few interventional experiments would be required to identify the true causal DAG. Finally, the last row suggests that the ratio of the number of MECs of size~1 to the total number of MECs also converges to $\approx 1/4$. Importantly, this would imply that $\approx 1/4$ of all causal DAGs can be uniquely identified without any interventional data. While a combinatorial analysis of the number of Markov equivalence classes for particular families of DAGs was initiated in~\citet{MEC_counting_1,MEC_counting_2}, these problems in general are wide open.

\begin{table}
\centering
	\begin{tabular}{ | c | c c c c c c c |}
	\hline	$p$ & 1 & 2 & 3 & 4 & 5 & 6 & 7 \\ \hline
	\# MEC	& 1 & 2 & 11 & 185 & 8782 & 1067825 & 312510571  \\
	(\# MEC)/(\# DAG)	& 1.00000 & 0.66667 & 0.44000 & 0.34070 & 0.29992 & 0.28238 & 0.27443 \\
	(\# MEC$_1$)/(\# MEC)	& 1.00000 & 0.50000 & 0.36364 & 0.31892 & 0.29788 & 0.28667 & 0.28068 \\ \hline	   
	\end{tabular}
	
	\begin{tabular}{ | c | c c c  |}
	\hline $p$ & 8 & 9 & 10 \\ \hline
\# MEC	& 212133402500 &  326266056291213 & 1118902054495975141 \\
(\# MEC)/(\# DAG)	&  0.27068 & 0.26888 & 0.26799\\
(\# MEC$_1$)/(\# MEC)	& 0.27754 & 0.27590 & 0.27507\\ \hline	   
\end{tabular}
\caption{The number of MECs, along with the  ratios of the numbers of MECs to DAGs and the ratios of the counts of MECs of size 1 (MEC$_1$) to the total number of MECs up to 10 nodes~\citep{Gillispie}.}
\label{table_MECs}
\end{table}

\subsection{Causal structure discovery algorithms and faithfulness}
\label{algorithms}

Since the overwhelming majority of available data  has been observational, most causal inference algorithms have been developed in this setting. A standard approach to causal structure discovery is \emph{constraint-based}, i.e., to treat causal inference as a constraint satisfaction problem with the constraints being the CI relations inferred from the data. A prominent example is the \emph{PC algorithm}, which starts in the complete undirected graph and iteratively removes edges $(i,j)$ if there exists $S\subset V\setminus\{i,j\}$ such that $X_i\independent X_j\mid X_S$. This results in the skeleton of the DAG; the immoralities are determined in a second step using the identified CI relations.

For such an algorithm to output the correct Markov equivalence class it is necessary that the inferred CI relations are \emph{faithful} to the true DAG. In particular, it has to hold that 
\begin{equation}
	\label{eq_faith}
	X_i\notindependent X_j\mid X_S \quad \textrm{for all } (i,j)\in E \textrm{ and all } S\subset V\setminus\{i,j\},
	\end{equation}
which is known as the \emph{adjacency faithfulness assumption}~\citep{adjacency_faithfulness}. Faithfulness violations can occur through cancellation of causal effects in the graph. Assumption (\ref{eq_faith}) seems harmless at first, since it is highly unlikely that causal effects in a DAG cancel each other out exactly. However, CI relations are inferred from data via hypothesis testing. So in the finite sample regime (\ref{eq_faith}) must be strengthened. In the Gaussian setting, where CI relations can be tested using partial correlations $\rho_{ij\mid S}$, (\ref{eq_faith}) leads to the definition of \emph{strong faithfulness}~\citep{strong_faithfulness}: $$\rho_{ij\mid S}\geq \lambda \quad \textrm{for all } (i,j)\in E \textrm{ and all } S\subset V\setminus\{i,j\},$$
where $\lambda\asymp \sqrt{\log(p)/n}$ to guarantee uniform consistency of the PC algorithm~\citep{Kalisch1}. 

	\begin{figure}[!t]
	\centering
	\includegraphics[width=4.0in,angle=0]{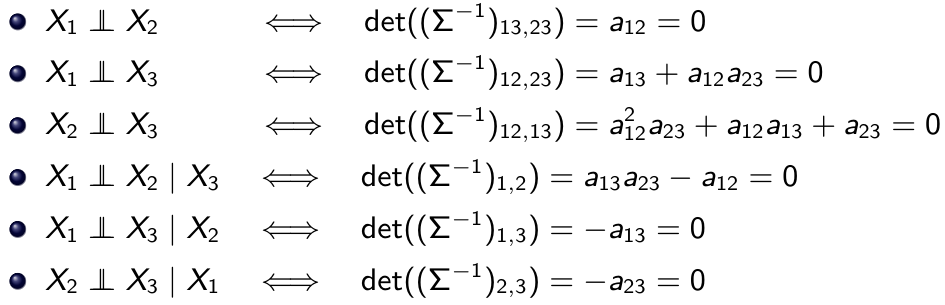}\qquad\includegraphics[width=1.5in,angle=0]{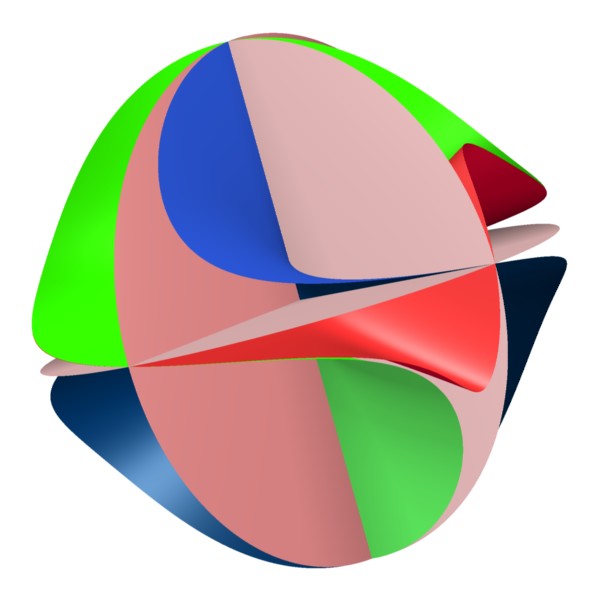}
	\caption{The unfaithful distributions for a 3-node fully connected DAG correspond to a collection of 6 hypersurfaces, each of which is defined by the vanishing of an almost principle minor; the 3 linear hypersurfaces are shown in pink and the 3 non-linear hypersurfaces in red, blue and green; illustration taken from~\citet{Buehlmann_Uhler}.}
	\label{fig_faithfulness}
\end{figure}

Since the strong faithfulness assumption is critical for the consistency of various prominent causal inference algorithms, it is important to understand how many samples are needed in general to satisfy it. Algebraic geometry has played a major role in answering this question~\citep{Buehlmann_Uhler,singular_faithfulness}. To see why, consider a Gaussian linear structural equation model on the fully connected DAG on 3 nodes with edges $1\to 2$, $1\to 3$, $2\to 3$. For simplicity, we assume that all error variances are equal to 1 and hence $(X_1,X_2,X_3)\sim\mathcal{N}(0,\Sigma)$, where $\Sigma^{-1}=(I-A)(I-A)^T$ and $A$ is strictly upper triangular containing the causal effects $a_{12}$, $a_{13}$ and $a_{23}$. Since no edge is missing, any CI relation is unfaithful to the DAG. On 3 nodes, there are 6 possible CI relations. For Gaussian distributions, any CI relation corresponds to the vanishing of an \emph{almost principle minor}, as shown in Figure~\ref{fig_faithfulness}. Hence, faithfulness violations correspond to a collection of \emph{real algebraic hypersurfaces} and understanding how restrictive the strong faithfulness assumption is requires the computation of the volume of \emph{tubes} around these hypersurfaces. This was achieved using tools from real algebraic geometry, namely Crofton's formula and Lojasiewicz inequality in~\citet{Buehlmann_Uhler} and using real log-canonical thresholds in \citet{singular_faithfulness}. These results were then used to compute the scaling of number of samples to number of variables that lead to the tubes filling up the whole space. This is important since in this case no faithful distribution exists. In the high-dimensional setting, this scaling was shown to be as bad as $p_n=o(\log n)$, a real limitation for the application of algorithms that rely on the faithfulness assumption, including the PC algorithm. These results also provide an example of how methods from algebraic geometry can be applied in the setting of high-dimensional statistics.



\subsection{DAG associahedra for causal inference from interventional data}
\label{DAG_associahedra}

With this understanding of unfaithful distributions as a collection of hypersurfaces, it is clear that obtaining algorithms with better consistency guarantees requires removing some of these hypersurfaces, i.e., testing fewer CI relations. Given a permutation (i.e., ordering) of the nodes $\pi$ that is consistent with the true DAG $G$ (i.e., if $i\to j$ in $G$, then $i< j$ in the ordering $\pi$), then by the Markov property $G$ can be recovered by testing only one CI relation per edge, namely the conditioning set consisting of all ancestors of $i$ and $j$ with respect to $\pi$, i.e.,
\begin{equation}
\label{eq_pi1}
X_i\independent X_j \mid X_S, \quad \textrm{ where }\, S=\{k\in V : k\leq i \textrm{ or } k\leq j \textrm{ w.r.t. } \pi\}\setminus\{i,j\}.
\end{equation}
The true ordering, however, is in general unknown and must be inferred from data. A natural approach following Occam's Razor is to associate to each permutation $\pi$ a DAG $G_{\pi}$ using (\ref{eq_pi1}) and to then return the \emph{sparsest permutation}, i.e., the sparsest DAG among all permutations. This approach is uniformly consistent under strictly weaker conditions than strong faithfulness, namely provided the sparsest DAG is in the true Markov equivalence class~\citep{SP_alg}. However, these improved consistency guarantees were achieved at a large computational price, since determining the sparsest permutation requires searching over all $p!$ permutations. 

	\begin{figure}[!t]
	\centering
	\includegraphics[width=6.0in,angle=0]{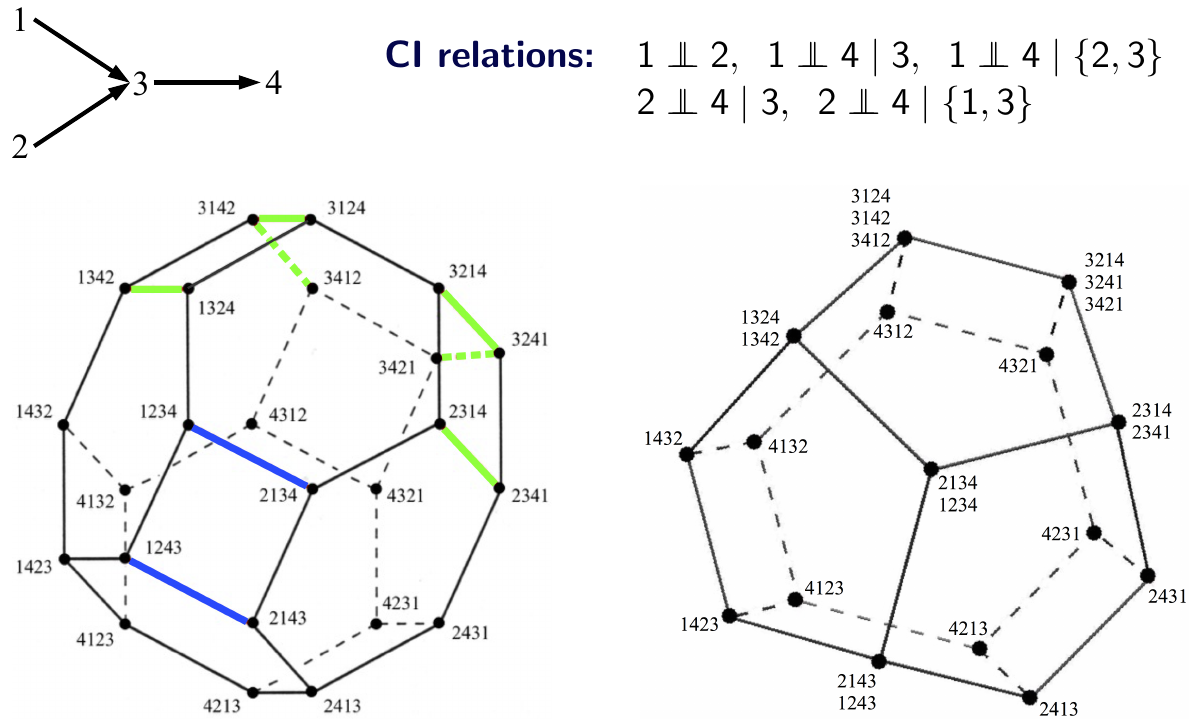}
	\caption{3-dimensional permutohedron consisting of all permutations of length 4 and the DAG associahedron for a particular 4-node DAG; illustration taken from~\citet{DAG_associahedra}.}
	\label{fig_DAG_ass}
\end{figure}

This raises the question whether replacing the exhaustive permutation search by a greedy search  could be used for causal inference. Greedy search algorithms are commonly applied for causal inference, most notably Greedy Equivalence Search, a greedy search over the space of Markov equivalence classes~\citep{Chickering}. The convex hull of all permutations of length~$p$ gives rise to a $(p-1)$-dimensional polytope, known as the \emph{permutohedron}, whose vertices are the permutations. Two permutations are connected by an edge in the permutohedron if and only if they differ by a neighboring transposition. The 3-dimensional permutohedron of all permutations of length 4 is shown in Figure~\ref{fig_DAG_ass}. It is shown in~\citet{GSP} that a greedy search in the permutohedron is consistent, i.e.~it outputs the correct Markov equivalence class when the sample size goes to infinity, under strictly weaker conditions than faithfulness. The sequences in Table~\ref{table_MECs} suggest that the number of MECs grows much faster than the number of permutations. Hence it is remarkable that a greedy search on the space of permutations has similar consistency guarantees as greedy search on the space of Markov equivalence classes, despite a large reduction in the search space.

In fact, the search space can be reduced further by identifying permutations whose DAGs $G_\pi$ and $G_{\pi'}$ are the same, since the number of edges in such graphs is necessarily the same. Such permutations are connected by edges in the permutohedron. Contracting these edges gives rise to a polytope~\citep{DAG_associahedra}, known as the \emph{DAG associahedron}, which can also be obtained by a different construction, namely by associating to each edge in the permutohedron a CI relation as described in~\citet{Morton_Pachter} and contracting all edges corresponding to CI relations in the underlying DAG~\citep{DAG_associahedra}; see Figure~\ref{fig_DAG_ass}. Thus DAG associahedra are a generalization of the prominent (undirected) graph associahedra~\citep{Carr_Devadoss}, that are obtained from the permutohedron by contracting all edges corresponding to separations or CI statements in an undirected graphical model. Hence the quest for a causal inference algorithm that is consistent under strictly weaker conditions than faithfulness and as a consequence achieves higher accuracies than previous algorithms in the high-dimensional setting led to the development of DAG associahedra and new results in convex geometry that are of independent interest. 

Recent years have witnessed a paradigm shift in the kinds of data that is being collected. In genomics, but also various other application areas, large-scale interventional datasets are being produced by deliberately altering some components in the system, such as genes. By further reducing the search directions, the greedy sparsest permutation search algorithm described above was extended to the first provably consistent algorithm for causal inference from a mix of observational and interventional data~\citep{IGSP,soft_IGSP}. For an application of these algorithms to learning gene regulatory networks see e.g.~\citet{IGSP,diff_DAG,soft_IGSP}.

\subsection{Open problems and related literature}
\label{applications}

While faithfulness  is well-understood from a geometric perspective, it is an open problem in algebraic geometry/combinatorics to understand the assumptions needed for consistency of the sparsest permutation algorithm. This is of great interest, since it is conjectured  that these are the weakest assumptions that guarantee consistency of any algorithm for learning the true Markov equivalence class~\citep{SP_alg}. Other polyhedral approaches for causal inference have been described~\citep{Cussens_2016,Jaakkola_2010} and it would be interesting to better understand how they relate to each other. So far we only considered causal inference when all variables are observed. However, for applications in the social sciences, latent variables are ubiquitous. The FCI algorithm and its variants generalize the PC algorithm to the latent setting~\citep{pc_algorithm}. It is an open problem to generalize greedy permutation search algorithms to the setting with latent variables. In addition, while CI relations are the only constraints that act on structural equation models in the fully-observed setting, in the latent setting there are additional constraints such as the \emph{Verma constraints}~\citep{Richardson_Evans}. While a full algebraic description of these constraints is not known, for linear Gaussian structural equation models a large subset has recently been characterized as nested determinants~\citep{Drton_Robeva}. In addition to these equality constraints, there are  inequality constraints. Describing these is very challenging, as demonstrated by the ongoing search for the semi-algebraic description of the set of matrices of fixed non-negative rank~\citep{Piotr_nnrank2,KRS2015}, which correspond to simple latent tree models in the discrete setting. Explicit knowledge of the defining equations and inequalities is crucial to answer questions of identifiability (e.g.~\citet{AMR2009}) or model selection (e.g.~\citet{Drton_Lin_Weihs_Zwiernik,Evans_model_selection}). Finally, we return to the beginnings of algebraic statistics on experimental design~\citep{Pistone2001} to point out a critical problem in the era of interventional data, namely to decide which interventions to perform in order to gain the most information about the underlying causal system.

\section{Phylogenetics}
\label{chp:marta}
This section treats a particular class of directed graphical models with latent variables, namely phylogenetic trees. Algebraic tools have been used since the end of the 20th century to address problems in phylogenetics \citep{felsenstein1978,hendy1989,evans1993,hendy1994}. In particular, \citet{lake1987} and \citet{cavender1987} opened the door to the development of phylogenetic reconstruction methods based on the polynomial equations implied by a particular evolutionary model and tree structure. 
Below, we survey this approach and then describe the major impact that algebraic statistics has had on phylogenetic reconstruction, model selection, and identifiability.

More detailed introductions to algebraic phylogenetics can be found in \citet{AllmanRhodeschapter4}, \citet{pachter2005}, \citet{sethbook}, \citet{piotrbook} or \citet[$\S$7,$\S$8]{steelbook}.

\subsection{Phylogenetic reconstruction}\label{subsec_phylo1}

\begin{figure}[!t]
\begin{center}
\includegraphics[scale=0.4]{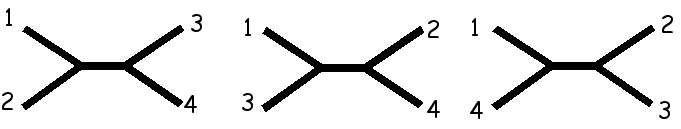}
\end{center}
\caption{\label{arbres4fulles} The three different unrooted tree topologies on four leaves are denoted as  $12|34$, $13|24$, and $14|23$ respectively.}
\end{figure}

A \emph{phylogenetic tree} is a tree graph whose leaves correspond to molecular sequences (for example, the whole genome of a species or a single gene) of different living species  and represents the speciation process that led to them: the interior nodes represent ancestral sequences, the edges evolutionary processes, and the leaves are labelled with the names of the living species. The \emph{topology} of a phylogenetic tree is the topology of the labeled graph; for example, Figure~\ref{arbres4fulles} shows the three different (unrooted) tree topologies for the species $\{1,2,3,4\}$, which are denoted as $12|34$, $13|24$, $14|23$.  While phylogenetic trees can be reconstructed using a variety of data including DNA or protein molecules, we here assume that the available data are a sequence  of characters $\texttt{A,C,G,T}$ (corresponding to the four nucleotides) of length $N$ for each living species in the tree.

In order to model the evolution of nucleotide data, it is convenient to assume that the substitution of nucleotides occurs randomly and following a Markov process on the phylogenetic tree $T$, where the internal nodes are latent variables. The state space of the random variables at the nodes of $T$ is $\{\texttt{A,C,G,T}\}$ and the parameters of the model are a distribution $\pi$ at a fixed interior node (which plays the role of the root) and the entries of the $4\times 4$ transition matrices $M^e$ associated to the edges $e$ of $T$. 
According to this hidden Markov process on $T$, where two nodes are independent given their least common ancestor, the probability of observing a character pattern at the leaves of the tree can be expressed as a polynomial in the model parameters. Assuming that the characters in each sequence have all evolved following the same evolutionary process and are independent of each other, the data are $N$ independent samples from a multinomial distribution.

Different restrictions imposed on the transition matrices give different evolutionary models: from the simplest  \emph{JC69} where $\pi$ is uniform and there is only one free parameter per edge (that is, on each edge, all conditional probabilities $P(x|y)$ are equal if $x\neq y$), to the  \emph{general Markov model} (\emph{GMM}) with no restrictions on the transition matrices or on $\pi$.

As an example, we consider the JC69 model on the tree $12|34$ shown in Figure \ref{arbres4fulles} relating the set of species $\{1,2,3,4\}$.
Denoting by $p_{x_1x_2x_3x_4}$ the joint probability of observing nucleotide $x_i$ at species $i$, it is not difficult to see (using the fact that the JC69 model is invariant under permutations of the four states) that
\begin{eqnarray}\label{modeleq}
p_{\tt AAAA} = p_{\tt CCCC} = p_{\tt GGGG} = p_{\tt TTTT},  \nonumber  \\
p_{\tt AAAC} = p_{\tt AAAG} = p_{\tt AAAT} = \dots = p_{\tt TTTG}  ,
\end{eqnarray}
whatever the parameters of the model. These are the first natural algebraic equations that arise from this type of evolutionary models. Although they are linear equations, they can be useful in model selection (see Section \ref{sec_modelsel}). They cannot be used to estimate the tree topology, because they hold for any joint distribution arising from a JC69 model on any of the trees in Figure \ref{arbres4fulles}; they are therefore called \emph{model invariants}.

Consider now
\begin{eqnarray}
p_{ACAC}+p_{ACGT}=p_{ACGC}+p_{ACAT} \quad \textrm{and} \label{lakeinvar}\\
p_{ACCA}+p_{ACTG}=p_{ACCG}+p_{ACTA}\label{lakeinvar2}.
\end{eqnarray}
Both \eqref{lakeinvar} and \eqref{lakeinvar2} hold for any set of JC69 parameters on the tree  $12|34$, but \eqref{lakeinvar} does not hold for all distributions on the tree $13|24$ and \eqref{lakeinvar2} does not hold for $14|23$ \citep{lake1987}. These equations that are satisfied for all joint distributions on a particular phylogenetic tree but not for all distributions on another tree are called \emph{topology invariants} \citep{steelbook}[8.3].

These equations were used by \citet{lake1987} to design a statistical test based on the $\chi^2$-statistic to infer the tree topology. These first attempts were not very successful \citep{huelsenbeck1995}, were only valid for simple models, and only used two of the relevant algebraic equations \citep{casfer2010}. The use of topology invariants for  phylogenetic reconstruction was thus halted, until the seminal work of \citet{allman2008}. We explain their main contribution in what follows.

Let $p \in \mathbb{R}^{256}$ be a distribution of character patterns on species $\{1,2,3,4\}$ as above, and consider its flattening matrix $flatt_{12\mid 34}(p)$ according to the split $12|34$, namely:

$ \qquad \qquad \qquad \qquad \qquad \qquad \qquad \qquad \qquad \scriptstyle{states \, \, at \, \, \, \, leaves \,\, 3 \, \, and \, \, 4}$ \vspace{-2mm}
$$flatt_{12\mid 34}(p) =\begin{array}{c}
\scriptstyle{states}\\
\scriptstyle{at}\\
\scriptstyle{leaves} \\
1,2\\
\end{array}
\left(
\begin{array}{ccccc}
p_{\tt AAAA} & p_{\tt AAAC} & p_{\tt AAAG} &  \dots & p_{\tt AATT} \\
p_{\tt ACAA} & p_{\tt ACAC} & p_{\tt ACAG} & \dots & p_{\tt ACTT}\\
p_{\tt AGAA} & p_{\tt AGAC} & p_{\tt AGAG} & \dots & p_{\tt AGTT} \\
\dots & \dots & \dots & \dots & \dots\\
p_{\tt TTAA} & p_{\tt TTAC} & p_{\tt TTAG} & \dots & p_{\tt TTTT}
\end{array}
\right).$$
\citet{allman2008} proved that if $p$ is a distribution from a Markov process on the tree $12|34$ (in the general Markov model), then $flatt_{12\mid 34}(p)$ has rank $\leq 4$; moreover it has rank $16$ if $p$ is a joint distribution on any of the two other trees $13|24$, $14|23$ (arising from a Markov process with generic parameters).

This result has allowed the development of new topology invariants under the general Markov model but, most importantly, the use of techniques such as rank approximation to propose methods to select the tree that best represents the data. This approach has been exploited in work that has attracted the attention of biologists \citep{SVDquartets, casfer2016}, and has allowed a generalization to the multispecies coalescent model and the use of these methods to estimate, not only gene trees, but also species trees \citep{chifmankubatko2015}. Some of these methods have been implemented in PAUP* \citep{paup}, one of the most widely used software packages in phylogenetics, which has opened the use of these tools to the biological community at large and has allowed the application in areas such as biodiversity preservation \citep{devitt2019}.

When dealing with real data, one only has access to a finite number of samples from the corresponding multinomial distribution.
Since many  methods are based on asymptotic tests, they may not be suitable for small samples.
Current research attempts to solve this issue and statistical tests based on algebraic tools are starting to be  developed for the finite-sample regime \citep{gaitherkubatko,sumner2017}.

These first approaches to algebraic phylogenetics have been restricted to trees on four species, but can be used in \emph{quartet-based methods} to infer large phylogenetic trees using only quartet data as input (that is, topologies of four species, in addition to an assessment score of the reliability of each quartet topology) \citep{strimmer1996,Ranwez2001,snirrao,davidson2015}. This approach has been used in \citet{casfer2016} to provide new support for the phylogenetic tree of eight species of yeast that was suggested by biological evidences and only obtained by certain reconstruction methods restricted to certain models. It would be of great interest to develop algebraic methods that can directly infer large trees; a first result in this direction is~\citet{sumner2017bis}.

The evolutionary models used in these algebraic approaches are more general than those commonly used  by biologists. Indeed, the usual approach in phylogenetics is to use a continuous-time Markov process. In this case, the transition matrix $M$ corresponding to an edge is of type $M=e^{tQ}$, where $Q$ is an instantaneous mutation rate matrix that operates for the duration $t\geq 0$. Not all transition matrices are of this type (i.e., not all Markov matrices are \emph{embeddable} in a continuous-time process); indeed, the logarithm of a transition matrix may not be real, and, if it is real, it may not be a rate matrix. \citet{RocaFernandez} proved that, for the Kimura 3-parameter model of nucleotide substitution, the  set of embeddable matrices represents only 9.4\% of all transition matrices. Moreover, it is commonly assumed that $Q$ is the same for all edges of the tree (i.e., the Markov process is \emph{homogeneous} in time), and that the process is stationary and time-reversible. This leads to one of the most used models in phylogenetics, the so-called \emph{general time-reversible} (\emph{GTR}) \emph{model}. While restricting to this model is quite controversial \citep{sumner2012} and might be too restrictive as we just insinuated, using GTR might be convenient because it considers less parameters than GMM (and hence the estimation of the parameters is more feasible). Algebraic approaches to phylogenetics avoid parameter inference altogether and make phylogenetic inference feasible for the most general Markov model, the GMM.

So far we have mainly focused on the recovery of the tree. As the number of trees grows super-exponentially in the number of leaves, accurately recovering the tree topology is a basic first step towards parameter recovery using methods such as maximum likelihood. Nevertheless, algebraic statistics can also lead to important results in obtaining estimates for continuous parameters of small phylogenetic trees. For example, using computational algebra one can compute the number of critical points of the log-likelihood function \citep{MLdegree} and then tools from numerical algebraic geometry can be used to obtain the global optimum \citep{KK}. Moreover, tools from computational algebra have provided major insights into the existence of a unique global optimum and provided analytical expressions to obtain it \citep{chorhendysnir2006, dinhmatsen}. In addition to equality constraints given by the model, the continuous parameters must also satisfy biological constraints and stochastic conditions, which are encoded as inequality constraints. While understanding these semi-algebraic constraints is difficult,  \citet{zwierniksmith}, \citet{matsen2009}, \citet{ARTaylor}, \citet{steelfaller} discuss which semi-algebraic constraints suffice to describe the model together with the algebraic constraints.

These algebraic tools are also being used in phylogenetic networks; for instance, \citet{chifmankubatko2019} introduce a new technique based on algebraic statistics to detect events in hybridization networks. These new tools for topology reconstruction are opening a new direction for phylogenetic reconstruction, with many interesting challenges from both statistical and algebraic points of view.

\subsection{Identifiability}

Although the use of algebraic statistics for proving identifiability of parameters of statistical models in phylogenetics is primarily of theoretical nature, it has been critical for proving the consistency of many phylogenetic reconstruction methods, including those based on likelihood. In the following, we provide a short overview.

\citet{chang1996} used algebraic tools  to prove that the GMM is \emph{generically identifiable}, that is, generic parameters are identifiable from the joint distributions of triplets of species (up to label swapping). The same holds for simpler models. Unfortunately, identifiability becomes much more involved for more complex models. 
Of particular interest are extensions 
that allow different sites to evolve at different rates, either by considering a $\Gamma$ distribution of rates across sites (but assuming that all sites evolve according to the same tree topology), or by considering a \emph{mixture model} (i.e., the joint distribution $p$ is a mixture of a certain number of distributions $p^i$ that have arisen from trees $T_i$ under a certain model $\mathcal{M}$, with unknown trees, continuous parameters, and mixing parameters). 

One well-known model that allows for different rates is the GTR+$\Gamma$, where  all sites evolve based on the same tree topology and with the same instantaneous mutation rate matrix, but the rate at which each site evolves follows a $\Gamma$ distribution (of certain fixed  parameters).
Although maximum likelihood estimation for this model has been widely-used by the biological community, identifiability of its parameters was established  using algebraic statistics only in 2008 in \citet{allmanGTRG}.

As far as mixture models are concerned, the first problem is to prove identifiability of the tree parameters. Without any constraints on the number of distributions, overfitting occurs and the continuous parameters are not identifiable. However, with appropriate constraints, the trees can be identified. Consider, for instance, mixtures on a single tree of four leaves evolving under the JC69 model. Equations \eqref{lakeinvar} and \eqref{lakeinvar2} are satisfied for all distributions on the tree $12|34$ and, as they are linear, they are also satisfied for any mixture of distributions on this tree. Therefore, as mixtures of distributions on any of the other two trees in Figure \ref{arbres4fulles} do not satisfy one of these equations, these topology invariants are able to identify the tree for this type of mixtures. In a similar way, the rank conditions mentioned in \ref{subsec_phylo1} allow a generalization that proves identifiability of trees in the case of mixtures on a single tree \citep{RhodesSullivant}. When mixtures on two different trees are considered, few positive results have been obtained \citep{allman2011identifiability}. For example, it is an open problem to determine whether, if one considers mixtures of GMM distributions on two trees $T_1$ and $T_2$,  the pair $\{T_1,T_2\}$ can be recovered from the mixed distribution.

Recently, algebraic tools have been used to prove the consistency of phylogenetic reconstruction methods for more complex models, including methods that reconstruct the species tree from gene trees according to the multispecies coalescent model. This is a very active area of research with important biological implications. See, for instance,  \citet{allmandegnanrhodes2018b}, \citet{allmanlongrhodes}, which all make use of deep algebraic tools.

Finally, in recent years there have been also incursions of algebraic statistics into questions about identifiability of phylogenetic and hybridization networks \citep{gross,chifmankubatko2019,banos}, but many open problems remain.

\subsection{Model selection}\label{sec_modelsel}

Another way in which algebraic statistics has been applied to phylogenetics is in selecting the evolutionary model that best fits the data. As mentioned above, a range of evolutionary models have been described ranging from JC69 to GMM (see for example the Felsenstein hierarchy in \citep{pachter2005}). Among them, the ones that have been deeply studied from an algebraic viewpoint are JC69, K80, K81, SSM and GMM. Following \citet{Kedzierska2012}, we explain here how the model invariants  for these models  can be used in model selection within the framework of phylogenetic mixtures.

\emph{Model invariants} on trees of $n$ leaves for an evolutionary model $\mathcal{M}$  are algebraic equations satisfied by all distributions arising from any set of parameters on the model $\mathcal{M}$ on any phylogenetic tree on $n$ leaves. The equations \eqref{modeleq} are an instance of model invariants for $\mathcal{M}$=JC69 on trees of four leaves. As they are linear equations, they are also satisfied for any mixture of a collection of distributions on these trees. In general, the space of mixtures of distributions on trees on $n$ leaves evolving under $\mathcal{M}$ (i.e. the set of mixtures of any number of distributions on trees on $n$ leaves evolving under $\mathcal{M}$)  is determined by the collection of linear model invariants. Moreover, the linear model invariants for $\mathcal{M}=$JC69, K80, K81, SSM are generated by binomial equations for any $n$ (analogous to \eqref{modeleq} and computed in \citet{casferked}), which leads to the exact computation of the likelihood maximum for data points coming from mixtures of distributions on a particular $\mathcal{M}$. Finally, these likelihoods can be combined into an information criterion for model selection (such as corrected Akaike or Bayesian Information Criterion). These tools were applied in \citet{Kedzierska2012} to real DNA data from the PANDIT database: while the usual model selection method chooses the most complex model GRT+$\Gamma$ (+invariable sites) and gives a tree incongruent with the accepted phylogeny, the method presented there selects a mixture of JC69 and leads to the accepted phylogenetic tree.

\section{Discussion}

Since its beginning in the late 1990s, the field of algebraic
statistics has grown rapidly.  The development of new theory and
algorithms for data analysis inspired by algebra, combinatorics and
algebraic geometry has brought together previously disconnected
communities of algebraists and statisticians. By now, algebraic
methods have touched on all major themes in statistics, such as
parameter identifiability and estimation, hypothesis testing, model
selection,
and Bayesian inference.  Conversely, problems and models from
statistics have inspired significant new ``pure'' developments in
algebraic combinatorics, high-dimensional commutative algebra, and
computational algebraic geometry. Various textbooks have been written on algebraic statistics: \citep{pachter2005,DSS09,sethbook,Aoki2012}, and for readers interested in using algebraic tools for statistical analysis, there is a package ``algstat"  implemented in R \citep{Ralgstat}.

We here provided an overview on developments made possible through the use of algebraic methods in three areas related to networks. We particularly focused on applications of algebraic statistics. However we did not touch upon many interesting developments of algebraic statistics. For example, significant contributions related to Markov bases have been applied to disclosure limitation~\citep{Fienberg2004} and genetics~\citep{MU_SNP}. Another
recent direction is the use of commutative algebra for experimental design in system
reliability~\citep{system_reliability}. Finally, another domain where algebraic techniques have been very fruitful is for the analysis of chemical reaction networks (see for example \citet{muller2016} and the work cited therein). It has been an exciting two decades for algebraic statistics. We have seen major impact of algebraic statistics on theoretical developments and, as summarized in this survey article, also on applications, and we expect that this discipline will expand into many further application domains.

\section*{ACKNOWLEDGMENTS}
MC was partially supported by AGAUR Project 2017 SGR-932 and MINECO/FEDER Projects MTM2015-69135 and MDM-2014-0445. CU was partially supported by NSF (DMS-1651995), ONR (N00014-17-1-2147 and N00014-18-1-2765), IBM, and a Sloan Fellowship. SP was partially supported by NSF (DMS-1522662) and IIT CISC (Center for Interdisciplinary Scientific Computation).

\bibliographystyle{natbib} 
\bibliography{algstats}

\end{document}